\documentclass{llncs}
\usepackage{makeidx}  
\usepackage{amsmath}  
\usepackage{bm}
\usepackage{latexcad}
\usepackage{url}
\begin{document}
\def\e{\mathrm e}                                               
\def\om{{\mathrm r}}                                            
\def\RR{\mathop{\cal R}\nolimits}                               
\def\cdc{,\ldots,}                                              
\def\1n{1\cdc n}                                                
\def\eq#1{\begin{equation}#1\end{equation}}                     
\def\eqs*#1{\begin{eqnarray*}#1\end{eqnarray*}}                 
\def\eqss#1{\begin{eqnarray}#1\end{eqnarray}}                   
\def\R{{\bbbr}}                                                
\def\T{{\xz\rm\scriptscriptstyle T}\xz}                         
\def\LW{{\rm\scriptscriptstyle LW}}                             
\def\s{{\hspace{.05em}\mathrm s}}                               
\def\xy{\hspace{.07em}}                                         
\def\xz{\hspace{-.07em}}                                        
\def\hy{\hspace{.04em}}                                         
\def\hz{\hspace{-.04em}}                                        
\def\ms{\mathstrut}                                             
\def\diag{\operatorname{diag}}                                  
\def\_#1{{^{}_{#1}}}                                            
\def\j{{\bar\jmath}}                                            
\def\i{{\bar\imath}}                                            
\def\jj{{\bar\jmath\bar\jmath}}                                 
\def\ii{{\bar\imath\hspace{.075ex}\bar\imath}}                  
\def\uv{{\hy\overline{\hz uv\hz}\hy}}                           
\def\vu{{\hy\overline{\hz vu\hz}\hy}}                           
\def\vv{{\hy\overline{\hz vv\hz}\hy}}                           
\def\u{{\hy\overline{\hz u\hz}\hy}}                             
\def\v{{\hy\overline{\hz v\hz}\hy}}                             
\def\l{\ell}                                                    
\def\dr{d^{\hspace*{.08em}\om}}                                 
\def\dLW{d^{\hspace{.05em}\LW}\xz}                              
\def\tdLW{\tilde d^{\hspace{.05em}\LW}\xz}                      
\def\Up#1{\vspace{-#1em}}                                       
\def\L{{\cal L}}                                                
\def\Ker{\operatorname{Ker}}                                    
\def\rank{\operatorname{rank}}                                  
\def\too{\!\to\!}                                               
\def\a{\alpha}
\def\x#1{{}}                                                    

\def\ovec{\overrightarrow}                                      
\def\G{G}                                                       
\def\PP{\mathop{\cal P}\nolimits}                               
\def\suml{\mathop{\sum}\limits}                                 
\def\maxl{\mathop{\max}\limits}                                 
\def\FF{\mathop{\cal F}\nolimits}                               
\def\FOto{\FF^{\bullet\to}}                                     
\def\FOij{\FF^{\,i\bullet\to j}}                                
\def\Fij{\FF_{ij}}                                              
\def\cd{\!\cdot}                                                

\pagestyle{headings}  
\addtocmark{Graph Metrics} 
%
%
%
\title{Studying new classes of graph metrics}
\titlerunning{New Classes of Graph Metrics}  
\author{Pavel Chebotarev
}
\authorrunning{Pavel Chebotarev}
%
\tocauthor{Pavel Chebotarev}
\institute{Institute of Control Sciences of Russian Academy of Sciences, 65~Profsoyuznaya~Str., Moscow 117997, Russia,\\
\email{pavel4e@gmail.com}}

\maketitle

\begin{abstract}
In data analysis, there is a strong demand for graph metrics that differ from the classical shortest path and resistance distances. Recently, several new classes of graph metrics have been proposed. This paper presents some of them featuring the cutpoint additive distances. These include the path distances, the reliability distance, the walk distan\-ces, and the logarithmic forest distances among others. We discuss a number of connections between these and other distances.
\end{abstract}

\section{Introduction}

For a long time, graph theorists studied only one metric for graph vertices, i.e., the shortest path metric~\cite{BuckleyHarary90}.
Gerald Subak-Sharpe was probably the first to investigate the resistance distance\footnote{In this paper, a \emph{distance\/} is assumed to satisfy the axioms of metric.} (electric metric)~\cite{Sharpe67a}. 
Later, this metric was studied in \cite{GvishianiGurvich87En} (see also \cite{Gurvich10DAM,Gurvich12RUTCOR}), \cite{ChandraRaghavanRuzzoSmolenskyTiwari89}, and~\cite{KleinRandic93}.

A distance $d(i,j)$ is called \emph{cutpoint additive\/} (or, \emph{graph-geodetic\/} \cite{Che11AAM}) whenever $d(i,j)+d(j,k)\!=\!d(i,k)$ if and only if $j$ is a cutpoint between $i$ and $j$ (i.e., every path connecting $i$ and $k$ visits~$j$).

In \cite{Che11AAM}, a general framework was presented for constructing cutpoint additive distances. Namely, it was shown that if a matrix $S=(s_{ij})$ produces a strictly positive \emph{transitional measure\/} on a graph $G$ (i.e., $s_{ij}\,s_{\!jk}\le s_{ik}\,s_{\!jj}$ for all vertices $i$, $j$, and $k,$ while $s_{ij}\,s_{\!jk}=s_{ik}\,s_{\!jj}$ if and only if every path from $i$ to $k$ visits~$j$), then the logarithmic transformation $h_{ij}=\ln s_{ij}$ and the inverse covariance mapping $d_{ij}=h_{ii}+h_{jj}-h_{ij}-h_{ji}$ convert $S$ into the matrix of a cutpoint additive distance. In the case of digraphs, five transitional measures were found\x{indicated} in \cite{Che11AAM}, namely,  the ``connection reliability'', the ``path accessibility'' with a sufficiently small parameter, the ``walk accessibility'', and two versions of the ``forest accessibility''.
The distances produced by the forest accessibility and walk accessibility on weighted multigraphs (networks) were studied in~\cite{Che11DAM} and~\cite{Che12DAM,CheDeza12,CheBapatBalaji12}, respectively.

\section{Notation}
\label{s_notat}

In the graph definitions we mainly follow~\cite{Harary69}.
Let $G$ be a weighted multigraph (a weighted graph where multiple edges are allowed) with vertex set $V(G)=V,$ $|V|=n>1,$
and edge set~$E(G)$. Loops are allowed; throughout the paper we assume that $G$ is connected. For brevity, we will call $G$ a \emph{graph}.
For ${i,j\in V,}$ let $n_{ij}\in\{0,1,\ldots\}$ be the number of edges incident to both $i$ and $j$ in~$G$; for every ${q\in\{\1n_{ij}\}}$, $w_{ij}^q>0$ is the weight of the $q\/$th edge of this type. Let
\eqs*{
a_{ij}=\sum_{q=1}^{n_{ij}}w_{ij}^q
}
(if $n_{ij}=0,$ we set $a_{ij}=0$) and $A=(a_{ij})_{n\times n};$ $A$~is the symmetric \emph{weighted adjacency matrix\/} of~$G$. 
All matrix entries are indexed by the vertices of~$G.$

The \emph{weighted Laplacian matrix\/} of $G$ is $L=\diag(A\bm1)-A,$ where $\bm1$ is the vector of $n$ ones and\x{} $\diag(\bm x)$ is the diagonal matrix with vector $\bm x$ on the diagonal.

\smallskip
For $v_0,v_m\in V,$ a $v_0\to v_m$ \emph{path\/} (\emph{simple path}) in $G$ is an alternating sequence of vertices and edges $v_0,\e_1,v_1\cdc\e_m,v_m$ where all vertices are distinct and each $\e_i$ is an edge incident to $v_{i-1}$ and~$v_i.$ The unique $v_0\to v_0$ path is the ``sequence''\,$v_0$ having no edges.

Similarly, a $v_0\to v_m$ \emph{walk\/} in $G$ is an \emph{arbitrary\/} alternating sequence of vertices and edges $v_0,\e_1,v_1\cdc\e_m,v_m$ where each $\e_i$ is incident to $v_{i-1}$ and~$v_i.$
The \emph{length\/} of a walk is the number $m$ of its edges (including loops and repeated edges). The \emph{weight\/} of a walk is the product of the $m$ weights of its edges. The weight of any subgraph $H$ of $G,$ denoted by $w(H),$ is defined similarly. The weight of a set of walks (of subgraphs) is the sum of the weights of its elements. By definition, for any vertex $v_0$, there is one $v_0\to v_0$ walk $v_0$ with length $0$ and weight~1.

Let $r_{ij}$ be the weight of the set $\RR^{ij}$ of all $i\to j$ walks in $G$, provided that this weight is finite. $R=R(G)=(r_{ij})_{n\times n}$ will be referred to as the \emph{matrix of the walk weights}.

\section{Cutpoint additive distances}
\label{s_cut-met}

By $d^\s(i,j)$ we denote the \emph{shortest path distance}, i.e., the length of a shortest path between $i$ and $j$ in~$G.$
This distance ignores the longer paths between $i$ and $j$ in $G$ as well as the numbers of paths of various lengths and edge weights.\x{}

A more sensitive graph metric is the resistance distance.
The \emph{resistance distance\/} between $i$ and $j$, $\dr(i,j),$ is the effective resistance between $i$ and $j$ in the resistive network whose line conductances equal the edge weights $w_{ij}^q$ in~$G$. The resistance distance has several expressions via the generalized inverse, minors, and inverses of the submatrices of the weighted Laplacian matrix of $G$ (see, e.g.,~\cite{KleinRandic93,Bapat99RD} or \cite{Sharpe67a,MooreSubak-Sharpe68} and the other papers by Subak-Sharpe and Styan published in the 60s and cited in~\cite{Che12DAM}). The resistance distance is proportional to the \emph{commute-time} (or, \emph{commute}) \emph{distance\/}~\cite{ChandraRaghavanRuzzoSmolenskyTiwari89}.

It turns out that these two types of graph distances are not sufficient for\x{} applications.
Say, in \cite{YenSaerensShimbo08}, for a parametric family of graph similarity measures (kernels) and a series of graph node clustering tasks, it was found that the results obtained for intermediate values of the family parameter $\theta$ are usually better than those obtained with the resistance/commute ($\theta\!\to\! 0$) and the shortest-path (${\theta\!\to\!\infty}$)\x{} kernels. Another example of a study where the resistance distance does not perform well is \cite{LuxburgRadlHein09}. The authors conclude: ``We prove that the commute distance converges [as the size of the underlying graph increases] to an expression that does not take into account the structure of the graph at all and that is completely meaningless as a distance function on the graph. Consequently, the use of the raw commute distance for machine learning purposes is strongly discouraged for large graphs and in high dimensions''.\x{}

Thus, for studying large networks, other sensitive graph metrics are required.

One of the new classes of such metrics are cutpoint additive distances.

\begin{definition}
\label{d_g-d}
For a multigraph $G$ with vertex set $V,$ a function $d\!:V\!\times\!V\to\R$ is called \emph{cutpoint additive\/} provided that $d(i,j)+d(j,k)=d(i,k)$ holds if and only if every path in $G$ connecting\/ $i$ and\/\x{} $k$ contains~$j$.
\end{definition}

Cutpoint additive functions can also be called \emph{bottleneck additive}; in \cite{Che11AAM,Che12DAM} they are referred to as \emph{graph-geodetic}.
A~theory of cutpoint additive distances was outlined in~\cite{Che11AAM}.

\begin{definition}
\label{def_trme}
Given a graph $G,$ a matrix\/ $S=(s_{ij})\!\in\!\R^{n\times n}$ is said to \emph{determine the transitional measure $s(i,j)\!=\!s_{ij},$ $i,j\!\in\! V,$ for $G$} if\/\x{} $S$
satisfies the \emph{transition inequality\/}\footnote{If $S$ has positive diagonal entries, then the transition inequality is equivalent to $s'_{ij}\,s'_{jk}\le s'_{ik},$ where $s'_{ij}
=\frac{s_{ij}}{\sqrt{s_{ii}\xy s_{\xz jj}}},\;$ $i,j,k\in V.$}
\eqs*{
s_{ij}\,s_{jk}\le s_{ik}\,s_{jj},\quad i,j,k\in V
}
and the \emph{graph bottleneck identity with respect to~$G\!:$}
\eqs*{
s_{ij}\,s_{jk}=s_{ik}\,s_{jj}
}
holds if and only if all paths in $G$ from $i$ to $k$ contain~$j$.
\end{definition}
The transition inequality is a multiplicative analogue of the \emph{triangle inequality for proximities\/} \cite{CheSha97,CheSha98} also called the ``unrooted correlation triangle inequality''~\cite{DezaLaurent97}.

\begin{theorem}[\cite{Che11AAM}]\label{t_rans}
Suppose that $S\!=\!(s_{ij})_{n\times n}$ determines a transitional measure for $G$ and has positive off-diagonal entries. Then $D=(d_{ij})_{n\times n}$ defined by
\eq{\label{e_rans}%
D=\tfrac{1}{2}(h{\bm1}^\T+\bm1 h^\T-H-H^\T),\; \mbox{where}\xy\footnote{This means that $h_{ij}=\ln s_{ij},$ $i,j=1\cdc n,$ i.e., the arrow indicates elementwise operations on matrices.}\:\:
H=\ovec{\ln S},
}
is a matrix of distances on\/~$V(G).$ Moreover$,$ this distance is cutpoint additive.
\end{theorem}

Consider two types of\x{} transitional measures in which for every $i\in V,\,$ $s_{ii}=1.$

\subsection{The path distances}

The \emph{path $\tau$-accessibility\/} of $j$ from $i$ in $\G$ is the total \emph{$\tau$-weight\/} of all paths from $i$ to~$j$:
\eq{
\label{e_path}
s_{ij}=w_\tau(\PP^{ij})=\sum_{P_{ij}\in\PP^{ij}}w_\tau(P_{ij}),
}
where $\PP^{ij}$ is the set of all $i\to j$ paths in~$\G,$
\eqs*{
w_\tau(P_{ij})=\tau^{l(P_{ij})}w(P_{ij}),
}
$l(P_{ij})$ and $w(P_{ij})$ are the length and the weight of $P_{ij},$ and $\tau>0.$

By definition, for every $i\in V,$ the unique ``path from $i$ to $i$'' is the path of length $0$ whose weight is unity, whence $s_{ii}=1,\;i=\1n$.

\begin{theorem}[\cite{Che11AAM}]\label{th_paths}
For any graph $\G,$ there exists $\tau_0>0$ such that for every\/ $\tau\in(0,\tau_0),$ $S=(s_{ij})$ defined by~\eqref{e_path}
determines a transitional measure for~$\G.$
\end{theorem}

\begin{corollary}[of Theorems~\ref{t_rans} and~\ref{th_paths}]\label{c_onpaths}
For any graph $\G,$ the matrix $D$ defined by \eqref{e_rans}$,$ where $S=(s_{ij})$ is defined by~\eqref{e_path} with a sufficiently small $\tau,$ determines a cutpoint additive distance on $V(G).$
\end{corollary}

\subsection{The reliability distance}
\label{s_conRel}

Consider a graph\footnote{In fact, the same results remain valid for strongly connected weighted digraphs.} $\G$ with edge weights $w_{ij}^p\in(0,1]$ interpreted as the intactness probabilities of the edges. Define $p_{ij}$ to be the $i\to j$ \emph{connection reliability\/}, i.e., the probability that at least one path from $i$ to $j$ remains intact, 
provided that the edge failures are independent. Let $P=(p_{ij})$ be the matrix of connection reliabilities for all pairs of vertices. For every $j\in V,\,$ $p_{jj}=1$, because the $j\to j$ path of length~$0$ is always intact.

The connection reliabilities can be represented as follows \cite[p.\,\x{}10]{Shier}:

\eq{
\label{e_Pways}
p_{ij}    =\suml_k        \!\Pr(P_k)
          -\suml_{k<t  }  \!\Pr(P_k P_t)
          +\suml_{k<t<l}\!\!\Pr(P_k P_t P_l)-\ldots
          +(-1)^{m-1}       \Pr(P_1 P_2\cdots P_m),
}
where $P_1, P_2\cdc P_m$ are all $i\to j$ paths in $\G$, $\Pr(P_k)=w(P_k),$ $\Pr(P_kP_t)=w(P_k\cup P_t)$, $P_k\cup P_t$ is the subgraph of $\G$ containing those edges that belong to $P_k$ or $P_t$, and so forth. By virtue of \eqref{e_Pways}, connection reliability is a modification of path accessibility that takes into account the degree of overlap for various paths between vertices.

\begin{theorem}[\cite{Che11AAM}]\label{th_reli}
For any graph $\G$ with edge weights $w_{ij}^p\in(0,1],$ the matrix $P=(p_{ij})$ of connection reliabilities determines a transitional measure for~$\G.$
\end{theorem}

\begin{corollary}[{of Theorems~\ref{t_rans} and~\ref{th_reli}}]
\label{co_paths}
For any graph $\G$ with $w_{ij}^p\in(0,1]$ for all edge weights$,$ the matrix of logarithmic distances $D$ defined by \eqref{e_rans}$,$ where $S=P=(p_{ij}),\x{}$ determines a cutpoint additive distance on~$V(G).$
\end{corollary}

For the following transitional measures, the diagonal elements $s(i,i)$ measure the (relative) strength of connections of every vertex to itself.

\subsection{The logarithmic forest distances}
\label{s_forest}

For a given $G$, by $\FF=\FF(G)$ and $\Fij=\Fij(G)$ 
we denote the set of all spanning rooted forests of $G$ and the set of all forests in $\FF$ that have vertex $i$ belonging to a tree rooted at~$j,$
respectively. Let
\vspace{-.2em}
\eq{
\label{e_fij}
f=w(\FF),\quad f_{ij}=w(\Fij),\quad i,j\in V(G);
}
by $F$ we denote the matrix $(f_{ij})_{n\times n}$; $F$~is called the \emph{matrix of forests\/} of~$G$.

Consider the matrix
$$
Q=(q_{ij})=(I+L)^{-1},
$$
where $I$ is the identity matrix.
By the matrix forest theorem \cite{CheSha97} (cf.\ Theorems 1--3 in \cite{ChungZhao10}), $Q$ does exist for any $G$\x{} and
\eq{\label{e_mft}
q_{ij}=\frac{f_{ij}}{f},\quad i,j=\1n.
}
Consequently, $F=fQ=f\cd\!(I+L)^{-1}$ holds.
$Q$~can be considered as a matrix providing a proximity measure for the vertices of~$G$\/ \cite{CheSha97,Che08DAM}.

\begin{theorem}[\cite{Che11AAM}]\label{th_fores}
For any $\G,$ the matrix $F=(f_{ij})$ defined by~\eqref{e_fij} determines a transitional measure for~$\G.$
Therefore$,$ the matrix $D$ defined by \eqref{e_rans} with $S=F$ determines a cutpoint additive distance on~$V(G).$
\end{theorem}

In \cite{Che11DAM}, a parametric family of logarithmic forest distances was studied. The family parameter determines a re-weighting of the edges of~$G.$ The marginal values of the parameter yield the shortest path and the resistance distances. The construction of the family is similar to that in the following section.

\subsection{The walk distances}
\label{s_prel}

For any $t>0,$ consider the graph $tG$ obtained from $G$ by multiplying all edge weights by~$t.$
If the matrix of the walk weights of $tG,\,$ $R_t=R(tG)=(r_{ij}(t))_{n\times n},$ exists, then\footnote{In the more general case of weighted \emph{digraphs}, the $ij$-entry of the matrix $R_t-I$ is called the \emph{Katz similarity\/}~\cite{Katz53} between vertices $i$ and~$j$.
}
\eqs*{
R_t=\sum_{k=0}^\infty(tA)^k=(I-tA)^{-1}.
}

Since $G$ is connected, $R_t$ is positive whenever it exists. The \emph{walk distances\/} are the distances appearing in the following theorem.
\begin{theorem}[\cite{Che11AAM}]\label{th_walk}
For any $\G,$ the matrix $R_t$ with $0\!<\!t\!<\!\rho^{-1},$ where $\rho$ is the spectral radius of $A$ determines a transitional measure for~$\G.$
Therefore$,$ $D$ defined by \eqref{e_rans} with $S=R_t$ determines a cutpoint additive distance on~$V(G).$
\end{theorem}

A topological interpretation of the walk distances in terms of the weights of walks and circuits in $G$ is given in~\cite{CheDeza12}.

\unitlength 0.85mm
\begin{figure}[ht!]
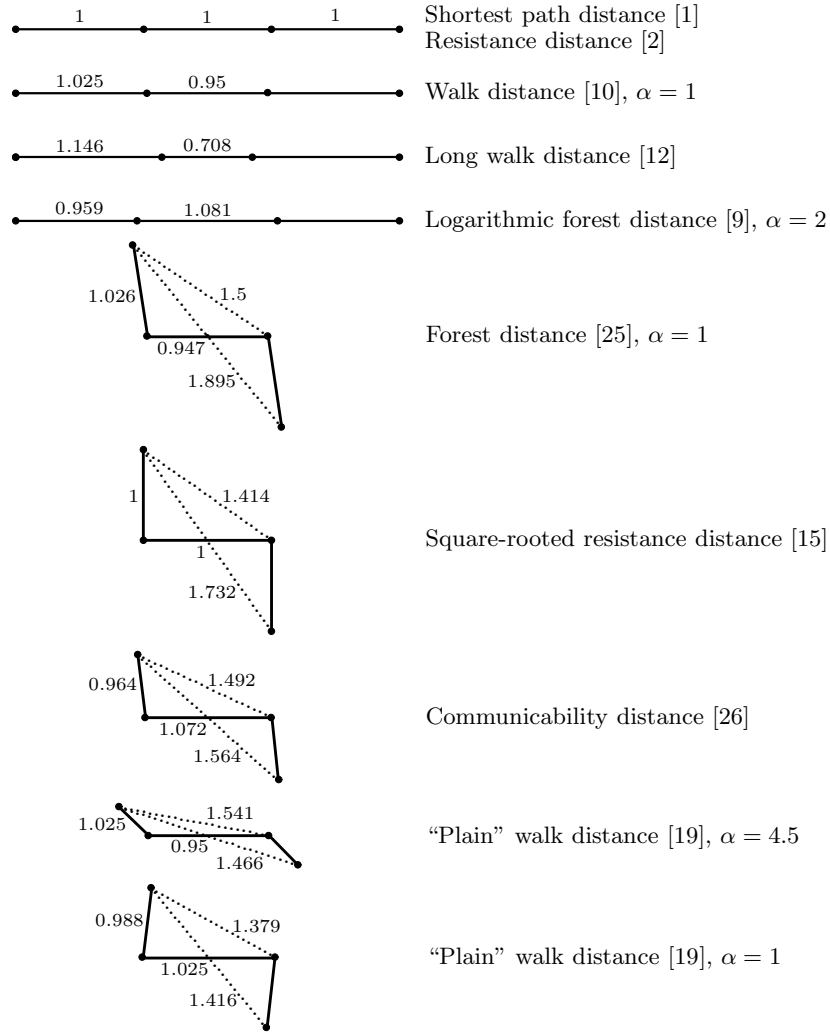

\input METRICS4.LP
\caption{\label{f_metrics}Various metrics on $P_4.$}
\end{figure}\subsection{The long walk distance}

The \emph{long walk distance\/} is defined \cite{Che12DAM} as the limit of the walk distance as $t\to(\rho^{-1})^-$ with a scaling factor depending on $t,$ $\rho,$ and~$n.$ Namely,
\eqs*{
\dLW(i,j)&\stackrel{\rm def}
 =&\lim_{t\to(\rho^{-1})^-}\frac{\ln r_{ii}(t)+\ln r_{jj}(t)-2\ln r_{ij}(t)}{n\rho^2(\rho^{-1}-t)},\quad i,j\in V.
}

This limit always exists and defines a cutpoint additive metric~\cite{Che12DAM}.
A number of closed-form expressions for $\dLW(i,j)$ are given in \cite{Che12DAM,CheBapatBalaji12}.
It turns out that the long walk distance is a counterpart of the resistance distance yielded by replacing the Laplacian matrix $L$ with the ``\emph{para-Laplacian matrix}'' $\L=\rho I-A.$
It is worth mentioning the rescaled version of the long walk distance introduced in \cite{CheBapatBalaji12}, namely, $\tdLW(i,j)\!\stackrel{\rm def}=\!n\xy\|p\|_2^2\,\dLW(i,j),$ where $p$ is the Perron vector of~$A.$\x{}

\section*{Conclusion}
In the paper, we discuss a number of connections between the metrics introduced above and other metrics.
Some applications of the cutpoint additive metrics in data analysis and machine learning are discussed in \cite{AlamgirVonLuxburg11a,KivimakiShimboSaerens12arX,SenelleGarcia-DiezMantrachShimboSaerensFouss13X}.\x{}
Fig.\,1 shows the values of several distances on the path graph $P_4$ with vertex set $\{1,2,3,4\}$ and edge set $\{(1,2),(2,3),(3,4)\};$\x{} the edge weights are unity. 
All distances are normalized in such a way that $d(1,2)+d(2,3)+d(3,4)=3;$ for all distances, $d(1,2)=d(3,4).$ The first four pictures correspond to cutpoint additive distances.
The remaining distances are listed in descending order of $d(1,4).$ Each picture presents the projection of the distance-obeying polygon onto the plane parallel to the line segments $(1,4)$ and $(2,3).$ It can be seen that some distances are shorter, other things being constant, 
for central vertices; some other distances are shorter for peripheral vertices.

Today, the task of finding a suitable graph metric for any particular application already looks quite realistic.\x{}


%
%
\bibliographystyle{splncs} 
\bibliography{all2}
\clearpage
\addtocmark[2]{Author Index} 
\renewcommand{\indexname}{Author Index}
\printindex
\clearpage
\addtocmark[2]{Subject Index} 
\markboth{Subject Index}{Subject Index}
\renewcommand{\indexname}{Subject Index}
\end{document}